\documentclass[preprint,twocolumn,5p,authoryear]{elsarticle}
  \def\appendixname{}
  \journal{Automatica}
\usepackage{graphicx}

\usepackage{mathtools,amsmath,amsthm,amsfonts,amssymb}
\usepackage{mathalfa}

\usepackage{epsfig} 
\usepackage{color}
\usepackage{bbm}
\usepackage{thmtools}
\usepackage{enumitem}
\usepackage{cancel}
\usepackage{psfrag}
\usepackage{pstool}
\usepackage{tikz,pgfplots}
\pgfplotsset{compat=1.17}

\usepackage{savesym}
\savesymbol{AND}
\usepackage{algorithm}
\usepackage{algorithmic}
\floatname{algorithm}{Algorithm}

\let\norm\undefined 
\DeclarePairedDelimiter\norm{\lVert}{\rVert}

\newcommand{\dd}{\textrm{d}}
\newcommand{\bracket}[1]{\left( #1 \right)}
\newcommand{\innprod}[2]{\ensuremath{\left\langle #1,#2 \right\rangle }}
\newcommand{\abs}[1]{\left| #1 \right|}
\newcommand{\BR}[1]{\ensuremath{\left\lbrace #1 \right\rbrace}}

\newcommand{\bbm}[1]{\left[\begin{matrix} #1 \end{matrix}\right]}
\newcommand{\sbm}[1]{\left[\begin{smallmatrix} #1
	\end{smallmatrix}\right]}
	
\renewcommand{\emph}[1]{{\it #1}}
	
\DeclareMathOperator*{\maxxx}{max}

\newcommand{\B}{\mathbb{B}}

\newcommand{\R}{\mathbb{R}}

\newcommand{\cA}{\mathcal{A}}

\newcommand{\cI}{\mathcal{I}}

\newcommand{\cK}{\mathcal{K}}

\newcommand{\cS}{\mathcal{S}}

\newcommand{\cU}{\mathcal{U}}
\newcommand{\cV}{\mathcal{V}}
\newcommand{\cW}{\mathcal{W}}

\newcommand{\cZ}{\mathcal{Z}}

\def\startmodif{\color{black}}
\def\stopmodif{\color{black}}

\theoremstyle{plain}
\newtheorem{lem}{Lemma}
\newtheorem{prop}{Proposition}
\newtheorem{cor}{Corollary}

\theoremstyle{definition}
\newtheorem{example}{Example}
\newtheorem{defn}{Definition}
\newtheorem{rem}{Remark}

\renewcommand\thmcontinues[1]{Continued}

\newcommand{\half}{\frac{1}{2}}

\DeclareMathOperator*{\argmin}{arg\,min}

\def\AT{\color{black}}
\begin{document}
\begin{frontmatter}

\title{Nearest Neighbor Control For Practical Stabilization of Passive Nonlinear Systems} 


\author[First,Second]{M. Z. Almuzakki}
\author[First]{B. Jayawardhana}
\author[Third]{A. Tanwani}

\address[First]{Engineering and Technology Institute Groningen, Faculty of Science and Engineering,\\ University of Groningen, The Netherlands (e-mail: {\tt\small \{m.z.almuzakki; b.jayawardhana\}@rug.nl}).}
\address[Second]{Department of Computer Science, Faculty of Science and Computer, Universitas Pertamina, Jakarta, Indonesia (e-mail: {\tt\small m.z.almuzakki@universitaspertamina.ac.id}).}
\address[Third]{Laboratory for Analyses and Architecture of Systems (LAAS) -- CNRS, Universit\'{e} de Toulouse, France (e-mail: {\tt\small aneel.tanwani@laas.fr}).}

\begin{abstract}                
This paper studies static output feedback stabilization of continuous-time (incrementally) passive nonlinear systems where the control actions can only be chosen from a discrete (and possibly finite) set of points. 
For this purpose, we are working under the assumption that the system under consideration is large-time norm observable and the convex hull of the realizable control actions contains the target constant input (which corresponds to the equilibrium point) in its interior. 
We propose a nearest-neighbor based static feedback mapping from the output space to the finite set of control actions, that is able to practically stabilize the closed-loop systems. 
Consequently, we show that for such systems with $m$-dimensional input space, it is sufficient to have $m+1$ \startmodif  discrete input points \stopmodif 
(other than zero for general passive systems or the target constant input for incrementally passive systems).
Furthermore, we present a constructive algorithm to design such $m+1$ \startmodif  nonzero \stopmodif input points that satisfy the conditions for practical stability using our proposed nearest-neighbor control.
\end{abstract}


\end{frontmatter}

\section{Introduction}\label{sec:intro}

{In several applications ranging from control of physical systems to networked control, exact implementation of a feedback control law is not possible due to the constraints at the level of sensors/actuators, or the constraints at the level of communication channels. Problems related to analysis, or the design of control laws, in the presence of such constraints have addressed considerable attention in the literature \citep{Persis2012,Delchamps1990,Elia2001,Hayakawa09,Jafarian2015a}. In this paper, we focus our attention on continuous-time dynamical systems where the input space is constrained to finite discrete sets.
}
An example of such a system with {\AT actuation} constraints is the design of the power take-off systems of the Ocean Grazer wave energy converter (WEC), where the device can only {\AT activate} a constant actuator systems from a pre-specified finite set \citep{Barradas-Berglind2016,Wei2017}.
\startmodif
Another example is the use of a fixed set of pulse attitude control thrusters for inducing 6 degree-of-freedom motion in spacecrafts or landers \citep{AretskinHariton2018,Kienitz2005}.
\stopmodif

Control design methods with appropriate analysis techniques, where binary input or minimal information is considered, have been discussed, among many others, in \citep{Elia2001,Kao2002} for linear systems, and in \citep{Cortes2006,Persis2012,Jafarian2015a} for the networked control systems setting. As these papers consider the use of binary input values per input dimension, the stabilization of an $m$-dimensional input-output system implies that there should be at least $2^m$ admissible input values and the stabilizing control law must dynamically assign one of these values as control input at every time {\AT instance}. 
In this paper, we shall focus on designing control laws with a minimal set of discrete control values whose cardinality is at most {\startmodif $m+1$, if we exclude the origin of the input space}.

We consider nonlinear systems described by
\begin{equation}\label{systems_eq}
\Sigma \ : \ \left\{
\begin{array}{rl}
\dot x         & = f(x) + g(x)u \\
    y         & = h(x) 
    \end{array}
\right.
\end{equation}
where the state $x(t)\in\R^n$ and the input and output signals $u(t),y(t)\in\R^m$.
The functions $f$, $g$, and $h$ are assumed to be continuously differentiable, $f(0)=0$, $g(x)$ is full-rank for all $x$, and $h(0)=0$. {For the developments carried out in this paper, the underlying assumption is that the input-output system $\Sigma$ is passive (in appropriate sense).
The basic problem studied in this paper is the stabilization of $\Sigma$ under limited actuation \slash information transmission; that is, the control input $u$ can only take values from a finite discrete set $\mathcal U:=\{u_0,u_1,u_2, \ldots , u_p\}$ with $u_i\in \R^m$ for each $i=0,\ldots,p$.}

Passive systems have received attention in different research fields as they are able to model physical phenomena exhibited by almost all thermo-chemo-electromechanical systems \citep{Ortega2013,Schaft2013}. In this regard, most of the aforementioned systems carry natural energy properties that can be related to passivity. In particular, such systems are said to be passive if the rate of change of the systems' ``stored energy'' never exceeds the power supplied by the environment through their external ports. There are different classes of passivity. For example, incremental passivity and differential passivity. These variations of passivity, along with the ``original'' passivity notion have been shown to be useful for control design purposes \citep{Jayawardhana2007,KOSARAJU2019466}. We refer interested readers to the various expositions on passive systems in \citep{Khalil2014,Ortega2013,Sepulchre2012,Schaft2016}. {\color{black} Our results are also applicable to a class of nonlinear systems that can be made passive by feedback control, as investigated, for instance, in \citep{Byrnes1991,Fradkov2008,Fradkov1998}. }

{\AT For the stabilization problem studied in this paper, it is assumed that we have a stabilizing output feedback law $y \mapsto F(y)$ (when $\mathcal U$ is continuum). When we impose the constraint that the actuation set $\mathcal U$ is finite,} two relevant questions for its stabilization are: a) how to map $F(y)$ to an element in $\mathcal U$?; and b) how to determine the minimal cardinality of $\mathcal U$? {\AT To address these questions for the system class $\Sigma$, we design a mapping $\phi:\R^{m} \to \mathcal U$, with $\mathcal U$  being discrete (and possibly minimal), such that $u=\phi(F(y))\in \mathcal U$ practically stabilizes $\Sigma$.} 

The question of designing the quantization mapping $\phi:\R^m \to \mathcal U$ has been addressed in various forms in literature. Since the input can only take the available values in the discrete set $\mathcal{U}$, the quantizer $\phi$, in some sense, defines the partition of the input space with respect to $\mathcal{U}$, {\AT where each cell of the partition is associated to an element of the set $\mathcal{U}$}. In most of the existing works, the input set $\mathcal{U}$ is chosen such that the resulting partition has some structure. {\AT For instance, when $\mathcal{U}:={\{-N,-N+1,\dots,N-1,N\}}^m$, a partition in the form of a regular grid facilitates design and analysis \citep{Ceragioli2007,Persis2012,Delchamps1990,Jafarian2015a,LibeHesp05,tatikonda2000control}. Other examples include logarithmic quantizers \citep{Elia2001,Fu2009}, which are optimal with respect to a certain density metric. However, if we fix the discrete set $\mathcal{U}$, then the question of finding the best possible partition for this given set $\mathcal{U}$ has not received much attention in the literature.}
In this paper, we address the later viewpoint by defining a simple static mapping that maps $F(y)$ to the nearest element in $\mathcal{U}$ which (practically) stabilizes the system. {\AT If certain stability conditions are satisfied, the resulting partition is described by convex polytopes, and to the best of authors' knowledge, such structures have not appeared in the literature on quantized control of nonlinear systems.}

The second question of finding the minimal set $\mathcal{U}$ for feedback stabilization has also received considerable attention. One question regarding this matter is on the minimal cardinality of the set $\mathcal U$. As an example, consider the work of \citep{NairFagn07}. In this paper, a discrete-time linear system, under some appropriate setting, is stabilizable if the number of bits per sample (rate of communication) is greater than the intrinsic entropy of the system. Similar results are available for continuous-time systems setting in \citep{Colo12,ColoKawa09}. {\AT To the best of authors' knowledge, there has not been a dedicated study on computing the entropy of passive nonlinear systems. Therefore, the question of how many symbols are necessary or sufficient for stabilization of a passive nonlinear systems has not been addressed. However, we do find some results on quantized control of passive system. In \citep{Cortes2006,Jafarian2015a}, under certain passivity structure in the dynamics, $\Sigma$ is shown to be practically stabilizable by using binary control for each input dimension which directly translates to $2^m+1$ elements in $\mathcal U$, e.g., $\mathcal U=\{0\}\cup \{-1,1\}^m$.}

As a relaxation of aforementioned results, and dealing with rather generic class of multi-input multi-output passive nonlinear systems, we show in this paper that such practical stabilization can be achieved by simply using $m+1$ elements in $\mathcal U$, in addition to $\{0\}$ or the required constant input $u^*$ when the system is required to track a desired constant reference $y^*$. {\AT We do so by proposing the nearest-neighbor based control laws and analyze the stability of the closed-loop systems when the input $u$ can only be taken from the finite discrete set $\mathcal{U}$. Moreover, we provide algorithmic procedure to construct minimal discrete sets that are able to practically stabilize the systems by means of nearest-neighbor based control law. 
Our design methodology is such that the overall closed-loop system is an interconnection of a passive system with an optimization-based selection rule for the input. Dynamical systems where the inputs are computed from solving an optimization problem, and are discontinuous appear in different applications \citep{BrogTanw20}. Passivity of the open-loop system is an important structural property that helps us analyzing the overall system in such cases.}
When quantization effect is of a particular concern, the interconnection of passive systems and quantizers has been studied for the past decade in various different contexts. For instance, the practical stability analysis of passive systems in a feedback loop with a quantizer using an adapted circle criterion for nonsmooth systems is presented in \citep{Jayawardhana2011}.

The rest of the paper is organized as follows. In Section~\ref{sec:prelims}, we provide some preliminaries on set-valued dynamics resulting from the use of nonsmooth control laws and on convex polytopes; and formulate the control problem. 
\startmodif
In Section~\ref{sec:result}, we describe our nearest neighbor control (NNC) approach, and the results showing practical convergence for passive systems. Using similar approach, we generalize the results in Section~\ref{sec:inc} by considering the practical stabilization of a nonzero equilibrium for incrementally passive nonlinear systems.
\stopmodif
Some simple designs of the minimal action set along with their construction procedures and properties associated to the NNC approach are provided and analyzed in Section~\ref{sec:minimal}. Finally, some concluding remarks are provided in Section~\ref{sec:conc}.

{\AT
A concise version of the results presented in Section~\ref{sec:result} has also appeared in the conference version of our paper~\citep{Jayawardhana2019nolcos}. However, in this article, we carry out the proofs differently and with more rigor, which allow us to tackle higher dimensional systems. The generalizations studied in Section~\ref{sec:inc}, and the design methods proposed in Section~\ref{sec:minimal} have not been addressed in any of authors' previous works.
}

\section{Preliminaries and Problem Formulation} \label{sec:prelims}

{\bf Notation:} For a vector in $\R^n$, or a matrix in $\R^{m\times n}$, we denote the Euclidean norm and the corresponding induced norm by $\| \cdot \|$. For a signal $z:\R_{\geq 0}\to \R^n$, the essential supremum norm of $z$ over an interval $I \subset \R_{\geq 0}$ is denoted by $\|z\|_I$. For any $c\in\R^n$, the set $\mathbb B_\epsilon(c) \subset \R^n$ is defined as, $\mathbb B_\epsilon(c) :=\{\xi\in \R^n | \|\xi-c\|\leq \epsilon\}$. For simplicity, we write $\mathbb B_\epsilon(0)$ as $\mathbb B_\epsilon$. The inner product of two vectors $\mu,\nu\in \R^m$ is denoted by $\langle \mu, \nu \rangle$. For a given set $\mathcal{S}\subset\R^m$, and a vector $\mu \in \R^m$, we let $\langle \mu, \cS \rangle := \{ \langle \mu, \nu \rangle \, \vert \, \nu \in \cS \}$. For a discrete set $\cU$, its cardinality is denoted by $\text{card}(\cU)$. The convex hull of vertices from a discrete set $\cU$ is denoted by $\text{conv}(\cU)$. The interior of a set $S\subset\R^n$ is denoted by $\text{int}\bracket{S}$. A unit vector whose $i$-th element is 1 and the other elements are 0 is denoted by $e_i$. A vector whose entries are 1 is denoted by $\mathbbm{1}$. A continuous function $\gamma:\R_{\geq 0}\to \R_{\geq 0}$ is of class $\cK$ if it is continuous, strictly increasing, and $\gamma(0)=0$. We say that $\gamma:\R_{\geq 0}\to \R_{\geq 0}$ is of class $\cK_\infty$ if $\gamma$ is of class $\cK$ and unbounded.

\subsection{Passive systems and observability notions}

The central object of this paper is the nonlinear control systems $\Sigma$ given in \eqref{systems_eq}.
The fundamental property that we associate with $\Sigma$ is that, it is {\it passive}, i.e., for all pairs of input and output signals $u,y,$ we have $\int_0^T \langle y(t),u(t) \rangle \dd t > -\infty$ for all $T > 0$; see \citep{Willems1972,Schaft2016,Ortega2013} for some primary references on passive systems. By the well-known Hill-Moylan conditions, the passivity of $\Sigma$ implies that there exists a positive definite storage function $H:\R^n\to\R_{\geq 0}$ such that $\langle \nabla H(x), f(x) \rangle \leq 0$ and $\langle \nabla H(x), g(x) \rangle = h^\top (x)$. Without loss of generality, we assume that the storage function $H$ is {\it proper}, i.e. all level sets of $H$ are compact. 

Using the passivity assumption on $\Sigma$, it is immediate to see that  $u\equiv 0$ implies that all level sets of $H$ are positively invariant. More precisely, for any $c>0$, if $H(x(0)) \leq c$ then $H(x(t)) \leq c$ for all $t\geq 0$. In other words, if we initialize the state of $\Sigma$ such that $x(0) \in \Omega_c:=\{\xi | H(\xi)\leq c\}$ with $u\equiv 0$ then $x(t)\in \Omega_c$ for all $t\geq 0$. We will use this property later to establish the practical stability of our closed-loop systems in conjunction with the following observability notion from \citep{Hespanha2005}. 
\startmodif
\begin{defn}
The system \eqref{systems_eq} is {\rm large-time initial-state norm observable} if there exist $\tau>0$, and $\gamma,\chi\in \mathcal K_\infty$ such that the solution $x$ of \eqref{systems_eq} satisfies
\[
\| x(t) \| \leq \gamma(\|y\|_{[t,t+\tau]})+\chi(\|u\|_{[t,t+\tau]})
\]
for all $t\geq 0$, $x(0)\in \R^n$, and locally essentially bounded and measurable inputs $u:\R_{\ge 0} \to \R^m$.
\end{defn}

In this work,
\stopmodif
we will use the large-time initial-state norm observability property for the autonomous system (with $u=0$):
\begin{equation}\label{autonomous_systems_eq}
       \dot x   = f(x), \quad
        \quad y  = h(x). 
\end{equation}
In this case, large-time initial-state norm observability of \eqref{autonomous_systems_eq} implies
\begin{multline}\label{eq:estAutxy}
    \exists \, \tau > 0,\gamma \in \mathcal K_\infty \ \text{such that, for each } x(0) \in \R^n,  \\
    \| x(t) \| \leq \gamma(\|y\|_{[t,t+\tau]}), \ \ \forall t\geq0.
\end{multline}

We note that in the standard passivity-based control literature, the notion of zero-state observability or zero-state detectability is typically assumed for establishing the convergence of the state to zero in the $\Omega$-limit set. However, these notions cannot be used to conclude the boundedness of the state trajectories given the bound on the output trajectories. Therefore, instead of using these notions, we will use the above large-time initial-state norm observability for deducing the practical stability based on the information on $y$ in the $\Omega$-limit set.

\begin{rem}\normalfont
If the dynamics in system~\eqref{autonomous_systems_eq} are linear, that is, $\dot x = Ax$, $y=Cx$, and the pair $(A,C)$ is observable, then one can quantify $\gamma$ in \eqref{eq:estAutxy} using the observability Gramian. In particular, if for $\tau > 0$
\[
W_\tau(t) = \int_{t}^{t+\tau} e^{A^\top (s-t)}C^\top Ce^{A (s-t)} \, \dd s
\]
then $x(t) = {\bracket{W_\tau(t)}}^{-1}\int_{t}^{t+\tau}e^{A^\top (s-t)}C^\top y(s) \, ds$, for each $t \ge 0$, and $\tau>0$, which in particular yields
\[
\norm{ x(t) }\le \norm{{\bracket{W_\tau(t)}}^{-1}} \int_{t}^{t+\tau} \!\! \norm{e^{A^\top (s-t)}C^\top} \, \dd s \!\!\! \sup_{s \in [t,t+\tau)}\vert y(s)\vert
\]
for each $t \ge 0$, and any $\tau > 0$.
\end{rem}

\subsection{Stabilization problem with limited control}
We are interested in feedback stabilization of the system $\Sigma$ described in \eqref{systems_eq} using the output measurements. The key element of our problem is that the input $u$ can only take values in a discrete set, which is finite. Thus, the objective is to find a reasonable way to map the outputs (taking values in $\R^m$) to a finite set such that the closed-loop system is stable in some appropriate sense. More formally, we address the following problem:

{\bf Practical output-feedback stabilization with limited control (POS-LC):} For a given system $\Sigma$ as in \eqref{systems_eq} and for a given ball $\mathbb B_\epsilon$ with $\epsilon > 0$, determine the {\AT finite set $\mathcal{U}:= \{u_0, u_1, \dots, u_p\} \subset \R^m$ with minimal cardinality}, and describe the mapping $\phi:\R^m\to \cU$ such that the closed-loop system of \eqref{systems_eq} with $u = \phi(y)$ satisfies $x(t) \to \mathbb B_\epsilon$ as $t\to\infty$ for all initial conditions $x(0) \in \R^n$.

In our problem formulation, both the construction of a discrete set $\mathcal U$, as well as the design of the stabilizing map $\phi$ constitute our control problem. Compared to the numerous works in the literature on quantized control, our job in solving POS-LC problem is facilitated under the passivity structure, along with the appropriate observability notion. In particular, for the first of results, we will work under the following basic assumption for solving POS-LC:
\begin{enumerate}
    \item[(A0)] The system $\Sigma$ in \eqref{systems_eq} is passive with a proper {\startmodif and positive definite} storage function $H$ and, the corresponding autonomous system~\eqref{autonomous_systems_eq} is large-time initial-state norm-observable for some $\tau >0$ and $\gamma \in \mathcal K_\infty$.
\end{enumerate}
{\AT
\begin{rem}\normalfont
In (A0), we require the storage function to be positive definite. In general, passivity of system~\eqref{systems_eq} only implies the existence of a positive semidefinite storage function. However, if we add zero-state-observability condition, then the resulting storage function is positive definite \cite[Lemma~1]{HillMoyl76}. In our setup, inequality \eqref{eq:estAutxy} implies such an observability notion.
\end{rem}
}
\subsection{Set-valued analysis: Basic notions}

In studying the aforementioned control problem, we recall some fundamental definitions found in the literature on differential inclusions and convex polytopes, which would be useful for analysis in later sections.

\subsubsection{Regularized differential inclusions}
It turns out that a mapping which maps output from a continuum to a discrete set of control actions is essentially discontinuous (with respect to usual topology on $\R^m$). Differential equations with such state-dependent discontinuities need regularization so that the solutions are properly defined. For a discontinuous map $F:\R^n\to\R^n$, we can define a set-valued map $\cK(F)$ by convexifying $F$ as follows
\[
\mathcal K(F(x)) := \bigcap_{\delta>0}\overline{\text{co}}(F(x+\mathbb B_\delta))
\]
where $\overline{\text{co}}(S)$ is the convex closure of $S$. The set-valued mapping $\cK(F)$ is the Krasovskii regularization of $F$, and under certain regularity assumptions on $F$, $\cK(F)$ is compact and convex-valued, and moreover it is upper semicontinuous.\footnote{A set-valued mapping $\Phi:\R^n \rightrightarrows \R^n$ is called {\it upper semicontinuous} at $x$ if for every open set $X$ containing $\Phi(x) \subset \R^n$, there exists an open set $\Xi$ containing $x$ such that for all $\xi\in \Xi$, $\Phi(\xi)\subset X$. Correspondingly, $\Phi$ is {\it upper semicontinuous} if it is upper semicontinuous at every point in $\R^n$.} For an upper semicontinuous mapping $\Phi:\R^n \rightrightarrows \R^n$, consider the differential inclusion
\begin{equation}\label{diff_incl}
\dot x \in \Phi(x)\qquad x(0)=x_0.
\end{equation}
A Krasovskii solution $x(\cdot)$ on an interval $I=[0,T),\ T>0$ is an absolutely continuous function $x: I \to \R^n$ such that  \eqref{diff_incl} holds almost everywhere on $I$. It is {\it maximal} if it has no right extension and it is a {\it global} solution if $I=\R_{\geq 0}$. For any upper semicontinuous set-valued map $\Phi$ such that $\Phi(\xi)$ is compact and convex for every $\xi\in\R^n$, the following properties have been established (see, e.g., \cite[Lemma 1]{Jayawardhana2011}): (i). the differential inclusion \eqref{diff_incl} has a solution on an interval $I$; (ii). every solution can be extended to a maximal one; and (iii). if the maximal solution is bounded then it is global. 

\subsubsection{Convex polytopes} 

Next, we present the definition of convex polytopes and some of their notable examples that are related to our problem. We refer to \citep{Okabe2009} and \citep{Toth2017} for additional material on this topic. In general, there are two basic representation of convex polytopes. Firstly, the vertex representation of a convex polytope in $\R^m$, or commonly referred to as the \emph{V-representation}, is an $m$-polytope defined by the convex hull of a finite set of points in $\R^m$; i.e.\ for any set of points $\mathcal{U}\subset\R^m$, the \emph{V-representation} of a convex polytope defined by $\mathcal{U}$ is given by $\mathcal{P}_{\rm{V}}(\mathcal U):={\rm conv}\bracket{\mathcal{U}}$. Another way to define an $m$-polytope is by intersecting finite-number of half-spaces, commonly referred to as the \emph{H-representation}, that is given by $\mathcal{P}_{\rm{H}}(A,b):=\BR{x\in\R^m | Ax\leq b }$. Note that both V-representation and H-representation of $m$-polytopes are equivalent, i.e.\ $\mathcal{P}_{\rm{V}}(\mathcal U)=\mathcal{P}_{\rm{H}}(A,b)$ with appropriate $A\in\R^{n\times m}$ and $b\in\R^n$. When it is clear from the context, we will omit the arguments in $\mathcal P_{\rm V}$ and $\mathcal P_{\rm H}$ in the rest of this paper.

One simple example of $m$-polytopes is the $m$-dimensional simplex, commonly referred to as $m$-simplex. For particular examples, 1-simplex is a line, 2-simplex is a triangle, and 3-simplex is a tetrahedron. The formal definition of $m$-simplices is given by:

\begin{defn}[$m$-simplex]
Let $\mathcal{S}:=\{s_0,s_1,\dots,s_m\}$ with $s_i\in\R^m$, $i=0,1,\dots,m$ be an affinely independent set, i.e.\ for any $s_i\in\mathcal{S}$, the set $\widetilde{\mathcal S}_i:=\{\tilde{s}\in\R^m \mid \tilde{s}=s_j-s_i,\forall s_j\in\mathcal{S}\setminus\{s_i\}\}$ is linearly independent. An $m$-simplex $\mathcal{S}_m$ is defined by,
\[
\mathcal{S}_m= {\rm conv}\left(\mathcal{S}\right):=\left\{ \sum\limits_{i=0}^{m} c_i s_i \biggm| \sum\limits_{i=0}^{m} c_i=1,\ c_i\geq 0 \right \},
\]
and we say that $b_{\cS_m} = \frac{1}{m+1} \sum_{i=0}^m s_i$ is its barycenter.
\end{defn}

\begin{example}\label{ex:m_simplex}
\normalfont
One special case of $m$-simplices is a \emph{regular} $m$-simplex $\mathcal{S}_{m,\rm{reg}}$ where all vertices have equal distances to its barycenter and, one possibly simple choice for such a simplex is 
\begin{equation}
    \mathcal{S}_{m,\rm{reg}}:={\rm{conv}}\bracket{\lambda \BR{ e_1,\dots,e_m,\frac{1-\sqrt{m+1}}{m}\mathbbm{1}}}\label{eq:reg_simplex}
\end{equation}
for some $\lambda \in\R_{>0}$.

Another notable example of $m$-polytopes is the $m$-dimensional hypercubes: the $m$-cubes and the $m$-cross-polytopes.
For a given $\lambda \in\R_{>0}$, an $m$-cube $\mathcal{C}_m$ is given by
$\mathcal{C}_m:=\BR{x\in\R^m\ |\ -\lambda \leq x_i\leq \lambda ;\ i=1,\dots,m},$
and an $m$-cross-polytope $\mathcal{C}_m^{\Delta}$ is given by
$\mathcal{C}_m^{\Delta}:=\rm{conv}\BR{\pm \lambda e_1,\dots,\pm \lambda e_m}.$
\end{example}

For our purposes, the utility of convex polytopes is seen in partitioning the output space $\R^m$ into a finite number of cells which can then be associated to a control action. In particular, given a finite set $\mathcal{S}\subset\R^m$ with ${\rm{card}}(\mathcal{S})=q$, the space $\R^m$ can be partitioned into $q$ number of cells where every cell contains all points in $\R^m$ that are closer to an element of $\mathcal{S}$ than any other element. Such cells are commonly referred to as Voronoi cells and are defined as follows.

\begin{defn}
Consider a countable set $\mathcal{S}\subset\R^m$. The Voronoi cell of a point $s\in\mathcal{S}$ is defined by
\[V_\mathcal{S}(s) := \BR{x\in\R^m\ |\ \|x-s\|\leq\|x-v\|,\ \forall v\in\mathcal{S}\setminus\{s\} }.\]
\end{defn}

\begin{rem}\normalfont
Note that every Voronoi cell is a closed and convex polyhedron since they can always be represented by the solution of a system of linear inequalities.
\end{rem}

\section{Nearest-Neighbor Control for Passive Systems}\label{sec:result}

In this section, we provide our first solution for the general passive systems when the practical stabilization of the origin is required. 
The motivation behind our design of these elements is to work with minimal number of elements in the set $\mathcal{U}$ which yield the desired performance using the static output feedback only. Toward this end, the only assumption we associate with the set $\mathcal U$ is the following:
\begin{enumerate}
    \item[(A1)] For a given set $\mathcal U:= \{u_0, u_1, u_2, \ldots, u_p\}$, with $u_0=0$, there exists an index set $\mathcal I\subset \{1,\ldots,p\}$ such that the set $\mathcal V:=\{u_i\}_{i\in \mathcal I}\subset \mathcal U$ defines the vertices of a convex polytope satisfying, $0\in\text{int}\bracket{\text{conv}\bracket{\mathcal{V}}}$.
\end{enumerate}

An immediate consequence of (A1) is the following lemma, which is used in the derivation of our forthcoming main result. 

\begin{lem}\label{lemma:1}
Consider a discrete set $\cU \subset \R^m$ that satisfies (A1). Then, there exists $\delta > 0$ such that
\begin{equation}\label{eq:bndVoronoi}
V_{\cU} (0) \subseteq \B_\delta,
\end{equation}
that is, the following implication holds for each $\eta \in \R^m$
\begin{equation}\label{eq:defCv}
    \| \eta\| > \delta \Rightarrow \ \exists \ u_i \in \cU \text{ s.t. } \|u_i + \eta \|< \| \eta\|.
\end{equation}
\end{lem}
\begin{proof}
Based on Assumption (A1), consider the sets
$\cI:=\{1,\dots, q\}$, 
\startmodif
and $\mathcal V:= \{v_1,\dots,v_q\} \subset \mathcal U$
\stopmodif
such that $q \le p$ and $0\in\text{int}\bracket{\text{conv}\bracket{\mathcal{V}}}$. Let $\mathcal{S}=\mathcal{V}\cup\BR{0}$. From the definition of Voronoi cells, it readily follows that $V_{\cU}(0) \subseteq V_{\cS}(0)$, and therefore, it suffices to show that $V_{\cS}(0) \subset \mathbb{B}_\delta$. Toward that end, we first observe that the Voronoi cell $V_{\cS}(0)$ can be described as
\begin{equation}\label{ineq:lemma1}
    V_\mathcal{S}(0):=\mathcal{P}_{\rm H}\bracket{\bbm{v_{1}&\dots&v_{q}}^\top,\ \frac{1}{2} \bbm{\norm{v_{1}}^2 &\dots&\norm{v_{q}}^2}^\top}.
\end{equation}
Thus, from~\eqref{ineq:lemma1}, we know that $V_\mathcal{S}(0)$ is a closed convex polyhedron. It remains to show that $V_\mathcal{S}(0)$ is bounded. Indeed, boundedness implies that we can choose $\delta = \maxxx\limits_{\Tilde{v}\in\ V_\mathcal{S}(0)}\bracket{\norm{\Tilde{v}}}$, such that $\mathbb B_{\delta}$ is the smallest ball containing the set $V_\mathcal{S}(0)$, which by definition of Voronoi cell is equivalent to \eqref{eq:defCv}.

To show that $V_\mathcal{S}(0)$ is bounded, we observe that, under (A1), there exists $\mu > 0$ such that $\mathbb{B}_\mu \subset {\rm conv}(\mathcal{V})$. Thus, for every $\tilde v \in V_\mathcal{S}(0)$, $\mu \frac{\tilde v}{\|\tilde v\|} \in {\rm conv}(\mathcal{V})$. Hence, there exist $\lambda_i \ge 0$ such that $\sum_{i=1}^q \lambda_i = 1$ and $\mu \frac{\tilde v}{\|\tilde v\|} = \sum_{i=1}^q \lambda_i v_i$. Consequently, from \eqref{ineq:lemma1}, it follows that
\[
\mu \frac{\tilde v^\top \tilde v}{\|\tilde v\|} = \sum_{i=1}^q \lambda_i v_i^\top \tilde v \le \frac{1}{2}\sum_{i=1}^q \lambda_i \|v_i\|^2
\]
and hence $\|\tilde v\| \le \frac{1}{2\mu}\sum_{i=1}^q \lambda_i \|v_i\|^2$.
\end{proof}

\begin{example}\label{example:exU}
\normalfont
A simple example of $\mathcal U$ in $\R^2$, satisfying (A1) is as follows:
\begin{equation}\label{ex:U}
\begin{aligned}
    \mathcal{U}_{\rm ex}&:=\alpha\BR{0,\sbm{\sin\bracket{\theta_{\rm ex}}\\\cos\bracket{\theta_{\rm ex}}},\sbm{\sin\bracket{\theta_{\rm ex}+\frac{2\pi}{3}}\\\cos\bracket{\theta_{\rm ex}+\frac{2\pi}{3}}},\sbm{\sin\bracket{\theta_{\rm ex}+\frac{4\pi}{3}}\\\cos\bracket{\theta_{\rm ex}+\frac{4\pi}{3}}}}\\
    &=:\BR{0,u_{{\rm ex},1},u_{{\rm ex},2},u_{{\rm ex},3}}
\end{aligned}
\end{equation}
with some $\theta_{\rm ex}\in\R$ and $\alpha\in(0,\infty)$. For this example, (A1) holds by taking $\mathcal{V}:=\mathcal{U}\setminus\{0\}$. Following the proof of Lemma~\ref{lemma:1}, we have $V_{\mathcal{U}}(0):={\rm{conv}}\bracket{\widetilde{\mathcal V}_{0}}$ where
\[\widetilde{\mathcal V}_{0}:=\alpha\BR{\sbm{\sin\bracket{\theta_{\rm ex}+\frac{\pi}{3}}\\\cos\bracket{\theta_{\rm ex}+\frac{\pi}{3}}},\ \sbm{\sin\bracket{\theta_{\rm ex}+\frac{3\pi}{3}}\\\cos\bracket{\theta_{\rm ex}+\frac{3\pi}{3}}},\ \sbm{\sin\bracket{\theta_{\rm ex}+\frac{5\pi}{3}}\\\cos\bracket{\theta_{\rm ex}+\frac{5\pi}{3}}}}.\]
Here, $\widetilde{\mathcal V}_{0}$ contains all vertices of the Voronoi cell $V_{\mathcal{U}}(0)$. Then, then the smallest $\delta$ that satisfies \eqref{eq:bndVoronoi} in Lemma~\ref{lemma:1} is given by $\delta=\alpha$. See Figure~\ref{fig:small_example} for an illustration.
\end{example}

\subsection{Unity output feedback}

\begin{figure}
    \centering
    \begin{tikzpicture}[scale=0.7, transform shape]
      \begin{axis}[
        xmin=-1.5, xmax=1.5,
        ymin=-1.5, ymax=1.5,
        ticks=none,
        width=1\linewidth,
        height=1\linewidth,
        color=black
        ]
        \addplot [color=blue!80!black,mark=*, ultra thick, mark size=3pt] coordinates {
		(0,0)} node[below left] {\large$0$} node[above=13,left=2] {\color{black}\large$V_{\mathcal U}(0)$};
        \addplot [domain=0:4*pi/3,samples=3,mark=*,color=black!20!blue,only marks,variable=t,
        ultra thick, mark size=3pt]({(sin(deg(t)))},{(cos(deg(t)))})
        node[pos=0,above]{\large $u_{\text{ex},1}$}
        node[pos=0.5,right=1]{\large $u_{\text{ex},2}$}
        node[pos=1,left=1]{\large $u_{\text{ex},3}$}
        node[pos=0,above, left=7]{\color{black}\large $V_{\mathcal U}(u_{\text{ex},1})$}
        node[pos=0.5,below=10]{\color{black}\large $V_{\mathcal U}(u_{\text{ex},2})$}
        node[pos=1,below=10]{\color{black}\large $V_{\mathcal U}(u_{\text{ex},3})$};
        
        \addplot [domain=pi/3:7*pi/3,samples=4,variable=t,
        ultra thick]({(sin(deg(t)))},{(cos(deg(t)))});
        
        \addplot [domain=0:2*pi,samples=100,variable=t,
        ultra thick, dash pattern=on 6pt off 2pt on 3pt off 2pt, color=red]({(cos(deg(t)))},{(sin(deg(t)))});
        \addplot [domain=1:1.2,samples=2,variable=t,stealth-,ultra thick]({(t*sin(deg(7.5*pi/8)))},{(t*cos(deg(7.5*pi/8)))}) node[pos=1,below] {\color{red}\large$\mathbb B_\delta$};
		\addplot [domain=0:1,samples=3,color=red!70!black,variable=t,
        ultra thick]({(t*sin(deg(pi/4)))},{(t*cos(deg(pi/4)))}) node[midway,below,sloped]{\large $\delta=\alpha$};
        
        \addplot [domain=1:5,samples=4,variable=t,
        ultra thick]({(t*sin(deg(pi/3)))},{(t*cos(deg(pi/3)))});
        \addplot [domain=1:5,samples=4,variable=t,
        ultra thick]({(t*sin(deg(3*pi/3)))},{(t*cos(3*deg(pi/3)))});
        \addplot [domain=1:5,samples=4,variable=t,
        ultra thick]({(t*sin(deg(5*pi/3)))},{(t*cos(5*deg(pi/3)))});
        \end{axis}
    \end{tikzpicture}
    \caption{\small Illustration of the nearest neighbor regions, or Voronoi cells, for discrete control actions set $\mathcal{U}_{\rm ex}=\BR{0,u_{\rm ex,1},u_{\rm ex, 2},u_{\rm ex,3}}$ given in Example~\ref{example:exU}. The triangular region around the origin is $V_{\cU}(0)$ and {\startmodif it is contained in $\mathbb{B}_\delta$, for some 
    $\delta=\alpha > 0$\stopmodif}.}
    \label{fig:small_example}
\end{figure}
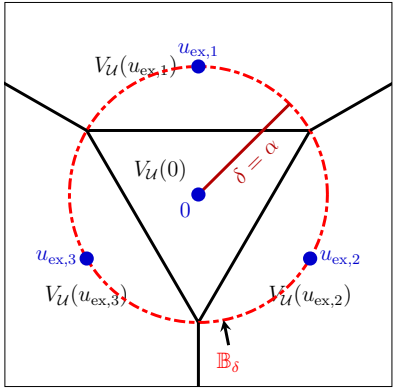

Using the result of Lemma~\ref{lemma:1} and the assumptions introduced thus far, we can define a feedback mapping $\phi$ which maps the measured outputs to the discrete set $\mathcal U$ to achieve practical stabilization. In this regard, we first consider the mapping $\phi:\R^m \rightrightarrows \cU$, defined as
\begin{equation}\label{phi_eq}
\phi(y) := \argmin_{v\in \mathcal{U}}\left\{\|v+y\|\right\}.
\end{equation}
The feedback control $u = \phi(y)$, with $\phi$ given in \eqref{phi_eq}, can be seen as a quantized version of the unity output feedback.
\startmodif
That is, when $\mathcal U$ is the continuum space $\R^m $, solution to the optimization problem~\eqref{phi_eq} is none other than $u=\phi(y)=-y$.
\stopmodif 
This quantization rule
\startmodif
$\phi$
\stopmodif
maps $-y$ to the nearest element in the set $\cU$ with respect to the Euclidean distance. The partitions in the output space induced by such quantization rule indeed result in Voronoi cells, and the resulting control law is hence discontinuous taking constant value in each of the Voronoi cells, see Fig.~\ref{fig:small_example}. By choosing $u = \phi(y)$, the closed system is thus given by
\begin{align}\label{nonsmooth_closedloop}
 \dot x & = f(x) + g(x)\phi(y) \\
    \nonumber y & = h(x).
\end{align}
As $\phi(y)$ is a non-smooth operator, we consider instead the following regularized differential inclusion 
\begin{align}\label{diff_incl_closedloop}
     \dot x & \in \mathcal K\big(f(x) + g(x)\phi(y)\big) = f(x) + g(x) \mathcal K(\phi(y)) \\
    \nonumber y & = h(x).
\end{align}
We note that the solution of \eqref{nonsmooth_closedloop} is basically interpreted in the sense of \eqref{diff_incl_closedloop}. In the following result, we analyze the asymptotic behavior of the solutions of \eqref{diff_incl_closedloop} and show that they converge to $\B_\epsilon$, for a given $\epsilon > 0$, if the elements of set $\cU$ satisfy certain conditions.

\begin{prop}\label{prop_nearest_control}
Consider a nonlinear system $\Sigma$ described by \eqref{systems_eq} that satisfies (A0), and a discrete set $\mathcal U \subset \R^m$ satisfying (A1) so that \eqref{eq:bndVoronoi} holds for some $\delta=\alpha>0$.
For a given $\epsilon > 0$, assume that 
\begin{equation}\label{eq:condEpsBall}
\gamma(\delta) \leq \epsilon.
\end{equation}
Then the control law $u=\phi(y)$, with $\phi$ given in \eqref{phi_eq}, globally practically stabilizes $\Sigma$ with respect to $\mathbb B_\epsilon$, that is, $\limsup_{t \to \infty} \vert x(t) \vert \le \epsilon$.
\end{prop}

\begin{proof}
For a fixed $y \in \R^m$, suppose that $\phi(y)=\{u_i\}_{i\in J_y}$ for some 
\startmodif
$J_y\subset \{0,1,\ldots, p\}$. 
\stopmodif
It follows from  \eqref{phi_eq} that $\{u_i\}_{i\in J_y}$ are the closest points to $-y$. Now, for each $i\in J_y$, we have that\footnote{\startmodif When $u_i$ is the closest point to $-y$, we know that the inequality $\|u_i+y\|^2\le \|u_j+y\|^2$ holds for all $j\in\{0,1,\ldots, p\}$. By taking $u_j=0$, and noting that $\|u_i+y\|^2 = \langle u_i+y,u_i+y\rangle = \|u_i\|^2 + 2\langle u_i,y\rangle + \|y\|^2$, we can also conclude that $\langle u_i,y \rangle \le -\frac{1}{2} \| u_i\|^2$.} 
(see \eqref{ineq:lemma1} also 
\startmodif
for construction of the following inequality using the definition of a half-space bounded by a hyperplane),
\stopmodif
\[
\langle u_i,-y \rangle \ge \frac{1}{2} \langle u_i,u_i\rangle \implies \langle u_i,y \rangle \le -\frac{1}{2} \langle u_i,u_i\rangle \le 0.
\]
Therefore, for each $y\in\R^m$, and $u_i \in \phi(y)$, $i \in J_y$, we get 
\begin{equation}\label{eq:bndUpLowInnProd}
-\|u_i\| \cdot \|y\| \le \langle u_i,y \rangle \le -\frac{1}{2} \|u_i\|^2
\end{equation}
Based on this property of $\langle \phi(y),y \rangle$, we can now analyze the behavior of the closed-loop system given by \eqref{diff_incl_closedloop}.

{\AT For the storage function $H$ associated with the open-loop system, we evaluate its derivative along the solutions of \eqref{diff_incl_closedloop} in following two cases:}\\
{\bf (i):} $ 0 \not \in \phi(y) = \{u_i\}_{i\in J_y}$ 
\startmodif
so that $J_y \subset \{1,\ldots, p\}$.
\stopmodif
Let $\cW_y := \phi(y)$, then
\begin{align*}
    \dot H(x) & = \langle \nabla H(x),\dot x \rangle \in \langle \nabla H(x) , f(x)+g(x)\mathcal K(\phi(y)) \rangle \\ 
    & = \langle \nabla H(x),f(x)\rangle + \langle y, \text{conv}(\cW_y) \rangle. 
\end{align*}
Based on the computation of $\langle \phi(y),y\rangle$, with non-zero $\phi(y)$, it follows that 
\begin{equation*}
\langle y, \text{conv}(\cW_y) \rangle \subset \big[- \|u_{y,\text{max}}\| \, \|y\|\, , \, - 0.5 \, \|u_{y,\text{min}}\|^2 \big],
\end{equation*}
where we let $\|u_{y,\text{max}}\|:=\text{max}_{w\in \mathcal W_y}\|w\|$, and  $\|u_{y,\text{min}}\|:=\text{min}_{w\in \mathcal W_y}\|w\|$.
Therefore,
\begin{align*}
  \dot H(x) \le -0.5 \, \|u_{y,\text{min}}\|^2; 
\end{align*}
when $0 \not\in \phi(y)$, or the other possibility is that,\\
{\bf (ii):} $ 0 = \phi(y)= \{u_0\} = \cW_y$ 
\startmodif
so that
$J_y = \{0\}$. 
\stopmodif
In this case, following the same arguments as in case {\bf (i)}
\[
\dot H(x) \in \langle \nabla H(x),f(x)\rangle + \langle y, \text{conv}(\cW_y) \rangle.
\]
Since $\{0\}$ is the only element of $\cW_y$, 
\begin{equation*}
\langle y, \text{conv}(\mathcal W_y) \rangle = \{0\}.
\end{equation*}
This implies that, for the case when $\phi(y) = \{0\}$, we have
\[
\dot H(x) = 0.
\]
{\AT Combining the two cases, it holds that for $J_y \subset \{0,1,\dots,p\}$, we have $\dot H(x) \le 0$, and $\dot H (x) = 0$, if and only if $0 \in \phi(y)$.} 
As $H(x)$ is non-increasing along system trajectories in both the cases {\bf (i)} and {\bf (ii)}, and since $H$ is proper, all system trajectories are bounded and contained in the compact set 
\startmodif
$\Omega_0 := \{ z \in \R^n\, | \, H(z) \le H(x(0))\}$. 
\stopmodif
Let 
\startmodif
$\cZ_x :=\{ z  \in \R^n \, \vert \, \phi(h(z)) = \{0\}\}$ 
\stopmodif
and let $M$ be the largest invariant set (with respect to system \eqref{diff_incl_closedloop}) contained in $Z_x$. By the LaSalle invariance principle, all trajectories belonging to the compact set $\Omega_0$ converge to the set $M$, see for example \cite[Theorem~6.5]{BrogTanw20}. 

We next show that, because of the large-time norm observability and Lemma~\ref{lemma:1}, it holds that $M \subset \mathbb{B}_\epsilon$. To see this, take an arbitrary point $z \in M$, and consider a solution of system~\eqref{diff_incl_closedloop} over an interval $[s,s+\tau]$ starting from $z$; that is, consider $x:[s,s+\tau] \to \R^n$ which solves \eqref{diff_incl_closedloop} and $x(s) = z \in M$. Due to the forward invariance of set $M$, the corresponding solution $x(t) \in M$, for each $t \in [s,s+\tau]$. Consequently,
\startmodif
$\phi(h(x(t))) = \{0\}$, 
\stopmodif
and because of Lemma~\ref{lemma:1}, $\vert h(x(t)) \vert \le \delta$ for each $t \in [s,s+\tau]$. Invoking the large-time initial state norm-observability assumption, it holds that $\| x(s) \| = \| z \| \le \gamma (\delta) \le \epsilon$, where the last inequality is a consequence of \eqref{eq:condEpsBall}. Since $z \in M$ is arbitrary, it holds that $M \subset \mathbb{B}_\epsilon$.

In summary, we have shown that
\[
x(t) \to M \subset \mathbb B_\epsilon \qquad \text{as} \ t\to\infty,
\]
for all initial conditions $x(0)\in\R^n$, and hence the desired assertion holds.
\end{proof}

As first application of Proposition~\ref{prop_nearest_control}, we are interested in specifying the invariant set when the set of control actions is described by a set of equidistant points along each axis of the output space.
\begin{cor}\label{cor_1}
Consider the system $\Sigma$ as in \eqref{systems_eq} satisfying (A0), and  $\mathcal U=\lambda\{-N,-N+1,\ldots,N-1,N\}^m$, with $\lambda > 0$ being the step size and $N$ a positive integer. Then the control law $u=\phi(y)$, where $\phi$ is as in \eqref{phi_eq}, globally practically stabilizes $\Sigma$ with respect to $\mathbb B_\epsilon$ where $\epsilon>0$ satisfies $\gamma(\lambda\sqrt{m})\leq \epsilon$.
\end{cor}

\begin{proof}
The proof follows {\it mutatis mutandis} the proof of Proposition \ref{prop_nearest_control}. The set $\mathcal U$ satisfies (A1) by taking 
\startmodif
$\mathcal V= \lambda\{-1,0,1\}^m\, \backslash \, \{0\}$. 
\stopmodif
It is also seen that $\delta = \lambda\sqrt{m}$, and by requiring $\gamma(\lambda\sqrt{m})\leq \epsilon$, all the hypotheses of Proposition~\ref{prop_nearest_control} hold.
\end{proof}

\begin{rem}
\normalfont
In contrast to the choice of $\mathcal{U}$ in Example~\ref{example:exU} where we used \eqref{ex:U} to construct the discrete set $\mathcal{U}$ in $\R^2$, the constant $\delta$ in Corollary \ref{cor_1} is less than $\maxxx_{\tilde v \in \widetilde{\mathcal V}} \| \tilde v\|$.
This is due to the choice of the set $\mathcal V$ in the proof of Corollary \ref{cor_1} that is dense enough such that $\{y\, | \, \phi(y)=0\} \subset \text{conv}(\mathcal V)$.  
From this corollary, one can conclude that two-level quantization with $N=1$ suffices to get a global practical stabilization property for passive nonlinear systems. This binary control law restricts however the convergence rate of the closed-loop system. It converges to the desired compact ball in a linear fashion and may not be desirable when the initial condition is very far from the origin. The use of higher quantization level (e.g., $N \gg 1$) can provide a better convergence rate when it is initialized within the quantization range.   
\end{rem}

\subsection{Sector bounded feedback}

We next present a generalization of the result in Proposition \ref{prop_nearest_control} on how the nearest neighbor rule can be used to quantize more generic nonlinear feedback laws. In Proposition~\ref{prop_nearest_control}, when $\mathcal U$ is the continuum space of $\R^m$, the resulting control law is simply given by $u=-y$, i.e., it is a unity output feedback law. Using standard result in passive systems theory, the closed-loop system will satisfy $\dot H \leq -\|y\|^2$.
\startmodif
Furthermore, the application of LaSalle invariance principle with zero-state detectability allows us to conclude that $x(t)\to 0$ asymptotically. 
\stopmodif
As the underlying system is passive, we can in fact stabilize it with any sector-bounded nonlinear feedback of the form $y \mapsto -F(y)$, where $F:\R^m \to \R^m$ satisfies
\begin{subequations}\label{eq:condSec}
\begin{gather}
k_1\|y\|^2\leq \langle F(y),y \rangle \leq k_2\|y\|^2, \quad 0<k_1\leq k_2 \\
\|F(y)\|\leq k_3\|y\|, \quad k_3 \ge k_1,
\end{gather}
\end{subequations}
for all $y \in \R^m$.
There are a number of reasons for considering such feedback laws rather than the unity output feedback law. For instance, we can attain a prescribed $L_2$-gain disturbance attenuation level or we can shape the transient behavior by adjusting the gains on different domain of $y$. 
In the following proposition, we consider such sector-bounded output feedback law $F(y)$, and how the nearest neighbor rule can be used to map such feedbacks in the limited control input set $\mathcal U$ to guarantee practical stabilization.

\begin{prop}\label{prop:phi_pass0_F}
Consider a nonlinear system $\Sigma$ described by \eqref{systems_eq} that satisfies (A0), and a discrete set $\mathcal U \subset \R^m$ satisfying (A1) so that \eqref{eq:bndVoronoi} holds for some $\delta>0$.
For the mapping $\phi$ given in \eqref{phi_eq}, {\color{black}let $\mu_{{\rm min},1}\in(0,1]$ be such that\footnote{{\color{black}The existence of such $\mu_{\text{min},1}$ is guaranteed by the assumption (A1) on $\mathcal U$.}}, for all $z\in\R^m$,
\begin{equation}\label{theta1max_eq}
\phi(z)\neq 0 \Rightarrow \langle \phi(z),-z \rangle \geq \|\phi(z)\|\|z\| \mu_{{\rm min},1}. 
\end{equation}
}
Assume that the constants $k_1, k_2, k_3$ describing the function $F$, as in \eqref{eq:condSec}, satisfy 
\begin{subequations}
\begin{gather}\label{sector_condition_F}
    {\color{black}\frac{k_1^2}{k_3^2} + \mu_{{\rm min},1}^2>1} \\
    \gamma\big(\delta/k_1\big) \leq \epsilon \label{eq:gammak_1}
    \end{gather}
\end{subequations}
for a given $\epsilon > 0$. 
Then the control law $u=\phi(F(y))$ globally practically stabilizes $\Sigma$ with respect to $\mathbb B_\epsilon$.
\end{prop}
\begin{proof}
We basically show that, for any $y\in \R^m$, we have
\begin{equation}\label{eq:secIneqMain}
\langle \phi(F(y)), y \rangle \in \{ -\kappa_{i,y}\|u_i\|\|y\| \, | \, i\in J_y\}
\end{equation}
for some 
\startmodif
$J_y \subset \{0,1,\ldots,p\}$
\stopmodif
such that $\phi(F(y)) = \{u_i\}_{i\in J_y}$ and $\kappa_{i,y} > 0$. The rest of the proof follows a pattern similar to that of Proposition \ref{prop_nearest_control}. 

{\AT First, with $\phi(F(y))=\{u_i\}_{i\in J_y}$, suppose that $0\notin \phi(F(y))$, so that $J_y \subset \{1,\dots,p\}$.} It follows from \eqref{phi_eq} that $\{u_i\}_{i\in J_y}$ are the closest points to $-F(y)$ which implies that
\begin{align}\label{mui1}
    \langle \phi(F(y)),-F(y) \rangle & \in \{ \|u_i\|\|F(y)\|\mu_{i,1} \, | \, i\in J_y\},
\end{align}
where $\mu_{i,1}>0$ is such that $\langle u_i,-F(y) \rangle = \|u_i\|\|F(y)\|\mu_{i,1}$. Under the given hypothesis, $\mu_{{\rm min},1} \leq \mu_{i,1}$ for each $i \in J_y$. On the other hand, 
\begin{align}\label{theta2max_eq}
    \langle -F(y), -y \rangle & = \|F(y)\|\|y\|\mu_{2}.
\end{align}
Since $k_1\|y\|^2 \leq \langle F(y),y \rangle$ and $\|F(y)\|\leq k_3\|y\|$, the minimum value of $\mu_{2}$ (with respect to all choices of $F$ that satisfy \eqref{eq:condSec}) is given by $\mu_{{\rm min},2}= k_1/k_3$.

Now, note that, in general, $\kappa_{i,y}\in [-1,1]$. It can be shown that if \eqref{theta1max_eq}, \eqref{sector_condition_F}, and \eqref{theta2max_eq} hold with $\mu_2\in [\mu_{\text{min,2}},1]$, then there exist $\kappa_{\text{min}}>0$ such that $\kappa_{i,y}\in [\kappa_{\text{min}},1]$. 
For each $y \in \R^m$ and $i \in J_y$, we introduce the Gram matrix $G_{i,y}$ as
\begin{align*}
    G_{i,y}=
    &\bbm{\innprod{-y}{-y} & \innprod{-y}{-F(y)} & \innprod{-y}{u_i} \\ 
\innprod{-y}{-F(y)} & \innprod{-F(y)}{-F(y)} & \innprod{-F(y)}{u_i}\\
\innprod{-y}{u_i} & \innprod{-F(y)}{u_i} & \innprod{u_i}{u_i}},
\end{align*}    
having the property that (see also \citep{Castano16}) $G_{i,y}\succcurlyeq 0$ and thus ${\rm det}(G_{i,y})\ge0$. This implies that
\begin{align*}
    0&\le {\|y\|}^2{\|F(y)\|}^2{\|u_i\|}^2
    +2\innprod{-y}{-F(y)}\innprod{-F(y)}{u_i}\innprod{-y}{u_i}\\
    &\quad-{\|y\|}^2{\innprod{-F(y)}{u_i}}^2-{\|F(y)\|}^2{\innprod{-y}{u_i}}^2\\
    &\quad-{\|u_i\|}^2{\innprod{-y}{-F(y)}}^2.
\end{align*}
By rewriting above inequality in terms of their respective norms and constants $\mu_{i,1},\mu_2,$ and $\kappa_{i,y}$, we have that, for each $u_i,\ i\in J_y$
\begin{align*}
    && \kappa_{i,y}^2-2\ \mu_{i,1}\ \mu_{2}\ \kappa_{i,y} &\le 1-(\mu_{i,1}^2+\mu_2^2)\\
    \Rightarrow&&
    \bracket{\kappa_{i,y}-\mu_{i,1}\ \mu_{2}}^2&\le 1-(\mu_{i,1}^2+\mu_2^2)+\mu_{i,1}^2\ \mu_2^2\\
    \Leftrightarrow&&
    \abs{\kappa_{i,y}-\mu_{i,1}\ \mu_{2}} &\le \sqrt{1-(\mu_{i,1}^2+\mu_2^2)+\mu_{i,1}^2\ \mu_2^2}.
\end{align*}
\startmodif
From the last inequality, we can prove whether $\kappa_{i,y}>0$ whenever condition \eqref{sector_condition_F} is satisfied, by only investigating the case where $\kappa_{i,y}\le\mu_{i,1}\ \mu_{2}$. The last inequality, paired with condition \eqref{sector_condition_F}, gives the following result
\stopmodif
\begin{align*}
    \kappa_{i,y}&\ge \mu_{i,1}\ \mu_{2}-\sqrt{1-(\mu_{i,1}^2+\mu_2^2)+\mu_{i,1}^2\ \mu_2^2}\\
    & = \mu_{i,1}\ \mu_{2}-\sqrt{(1-\mu_{i,1}^2)(1-\mu_2^2)}\\
    &\ge \mu_{\text{min},1}\frac{k_1}{k_3}-\sqrt{(1-\mu_{\text{min},1}^2)\left(1-\frac{k_1^2}{k_3^2}\right)}\\
    &= \mu_{{\rm min},1}\frac{k_1}{k_3} - \sqrt{1-\bracket{\mu_{{\rm min},1}^2+\frac{k_1^2}{k_3^2}}+\mu_{{\rm min},1}^2\frac{k_1^2}{k_3^2}}\\
    &>\mu_{{\rm min},1}\frac{k_1}{k_3}- \sqrt{\mu_{{\rm min},1}^2\frac{k_1^2}{k_3^2}}=0.
\end{align*}
Note that the above arguments hold for all $i\in J_y$, and \eqref{eq:secIneqMain} holds for some $\kappa_{i,y} > 0$.

{\AT Secondly, in case, $J_y = \{0\}$, we have $\phi(F(y)) = \{0\}$ and $\langle \phi(F(y)), y\rangle = \{0\}$. Thus, \eqref{eq:secIneqMain} holds trivially since $u_0 = 0$.}

{\AT Combining the two cases, we see that \eqref{eq:secIneqMain} holds for $J_y \subset \{0,1,\dots,p\}$.} Following the same line of arguments as in the proof of Proposition \ref{prop_nearest_control}, \eqref{eq:secIneqMain} implies that the storage function is nondecreasing along the solutions of the closed-loop system and the solutions converge to a set $M$, where $M$ is the largest invariant set contained in 
\startmodif
$\cZ_x :=\{ z  \in \R^n \, \vert \, \phi(F(h(z))) = \{0\}\}$.
\stopmodif
Hence for any trajectory starting with initial condition $x(s) = z \in M$, it holds that the corresponding output satisfies $\|F(y(t))\| \leq \delta$ for all $t \ge s$. Since $k_1\|v\|^2\leq \langle F(v),v \rangle \leq \|F(v)\|\|v\|$ holds for all $v\in \R^m$, it follows that $\|y(t)\| \leq \frac{\delta}{k_1}$ for all $t \ge s$.
By the property of large-time initial-state norm-observability of \eqref{autonomous_systems_eq}, it holds that,
\begin{equation*}
\|z\| = \|x(s)\|  \leq \gamma(\delta/k_1) \leq \epsilon \qquad \forall t \ge s
\end{equation*}
and this holds for each $z \in M$. Hence, $M \subseteq \mathbb{B}_\epsilon$ and in particular, each trajectory converges to $\mathbb{B}_\epsilon$ as $t \to \infty$.
\end{proof}

\begin{rem}
\normalfont
The condition \eqref{sector_condition_F} requires that the nonlinearity should lie in a relatively thin sector bound. When $F(y) = ky$, i.e, it is a proportional controller with a scalar gain $k>0$, then the condition \eqref{sector_condition_F} holds trivially{\color{black}, since $\mu_{{\rm min},1}>0$ and $\frac{k_1}{k_3}=\frac{k}{k}=1$}. Consequently, it follows from this proposition that we can make the practical stabilization ball arbitrary small by assigning a large gain $k$.  
\end{rem}

\subsection{An illustrative example}

\begin{figure}
    \centering
    \includegraphics[width=0.95\linewidth]{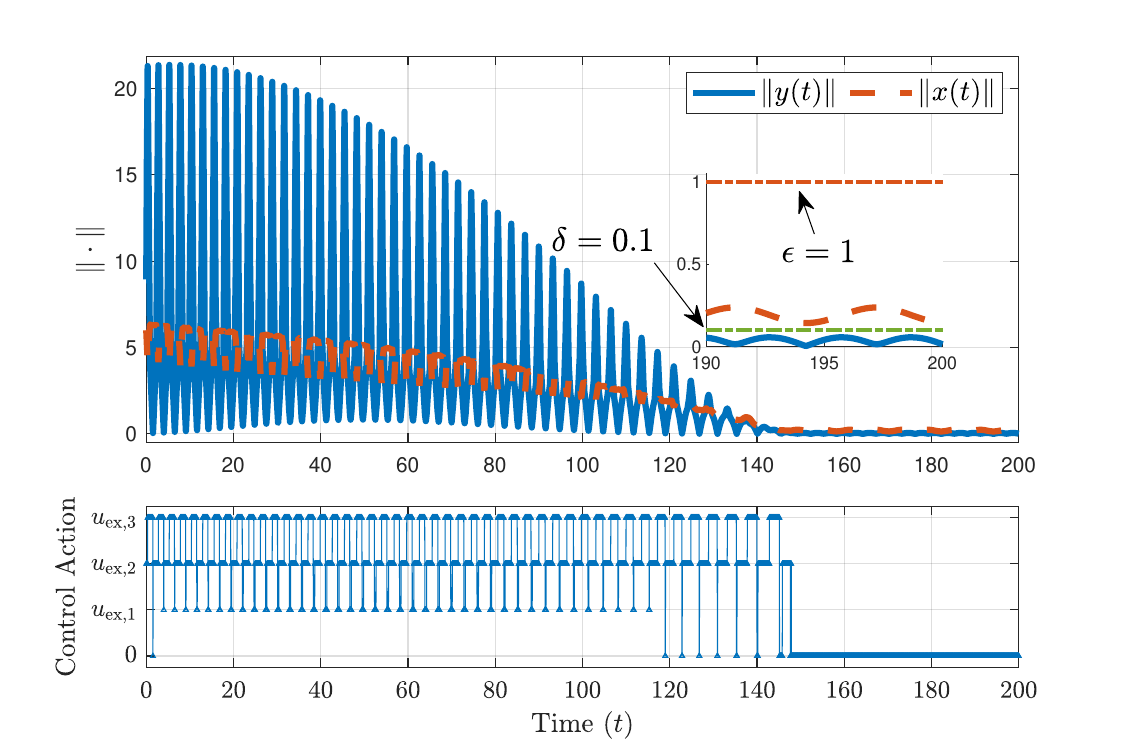}
    \caption{\small Simulation results of $\Sigma_{ex}$ using the control approach proposed in the Proposition~\ref{prop_nearest_control} with discrete input set $\mathcal{U}_{ex}$ as in \eqref{ex:U} and fixed parameters $\theta_{\rm ex}=0$ and $\alpha=0.1$. It can be seen that once both the state $x$ and the output $y$ enters their respective convergence ball, the control input is zero.}
    \label{fig:ex1}
\end{figure}

\begin{example}\label{ex:1}
\normalfont
Consider the following nonlinear system 
\begin{equation}
\Sigma_{\text{ex}} \ : \ \left\{    \begin{array}{rl}
\dot x         & = \bbm{-x_2  + x_3^3\\ x_1 \\ -x_1} + \bbm{1 & 0\\0&0\\0&1}u \\
y         & = \bbm{x_1\\x_3^3}
    \end{array} \right.
\end{equation}
where $x:=\bbm{x_1&x_2&x_3}^\top\in\R^{3}$ and $y:=\bbm{y_1&y_2}^\top,u:=\bbm{u_1&u_2}^\top\in\R^{2}$. It can be checked that by using the 
proper storage function $H(x)=\half x_1^2 + \half x_2^2 + \frac{1}{4} x_3^4$,
the system $\Sigma_{\text{ex}}$ is passive. Indeed, a straightforward computation gives us $\dot H = \langle y,u\rangle$. 
{\color{black} Note that the above system can be written as a nonlinear port-Hamiltonian system, describing a nonlinear RLC circuit \citep{Castanos2009}: $\dot x = J\nabla H(x) + gu$, $y=g^T\nabla H(x)$ where $J=\sbm{0 & -1 & 1\\1 & 0 & 0\\-1 & 0 & 0}$ and $g=\sbm{1 & 0\\0 & 0\\0 & 1}$. }  

We will now show that $\Sigma_\text{ex}$ satisfies the large-time initial-state norm observability condition. As the bound on $x_3$ for the large-time norm observability can directly be obtained from the output $y$, we need to compute the bound on $\sbm{x_1\\x_2}$. If we consider the sub-system of $\sbm{x_1\\x_2}$ with $x_1$ as its output (and is equal to $y_1$), it is a linear system with $A=\sbm{0 & -1\\1 & 0}$, $B=\sbm{1\\0}$, $C=\bbm{1&0}$ and its input is $x_3^3 = y_2$. Thus as $(A,C)$ is observable, the observability Gramian is given by
\[
W_{\pi}(t) = \int_{t}^{t+\pi} e^{A^\top (s-t)}C^\top Ce^{A(s-t)}\dd s = \frac{\pi}{2}\bbm{1 & 0\\0 & 1},
\]
whose inverse is simply given by $W_\pi^{-1}=\frac{2}{\pi}I_2$ and $\|W_\pi^{-1}\| = \frac{2}{\pi}$. Then for any $t>0$
\begin{equation*}
\sbm{x_1(t)\\x_2(t)} \!\!= \\ W_\pi^{-1}\int_{t}^{t+\pi}\!\!\!\!e^{A^\top (s-t)}C^\top \left(x_1(s)-\left(\! H*\sbm{x_3^3\\0}\right)\!(s)\!\right)\dd s,
\end{equation*}
where $*$ denotes the convolution operation and $H$ is the convolution matrix kernel given by $H(t)=Ce^{At}$. Since $\|e^{At}\|=1$ for all $t$, it follows then that
\begin{align*}
    & \left\|\bbm{x_1(t)\\x_2(t)}\right\| \\ & \leq  \|W_\pi^{-1}\|\int_{t}^{t+\pi} \!\! \left\|e^{A^\top (s-t)}C^\top \right\|\left\|x_1(s)-\left(H*\sbm{x_3^3\\0}\right)(s)\right\|\dd s \\
    & \leq \frac{2}{\pi} \pi \left(\|y_1\|_{[t,t+\pi]} + \|y_2\|_{[t,t+\pi]} \right) 
    \leq 4 \|y\|_{[t,t+\pi]}.
\end{align*}
Since by the definition of $y$, 
$\|x_3\|_{[t,t+\pi]}=\|y_2\|_{[t,t+\pi]}^{\frac{1}{3}}\leq \|y\|_{[t,t+\pi]}^{\frac{1}{3}}$, it follows from the inequality above that
\[
\|x(t)\| \leq 4 \|y\|_{[t,t+\pi]} + \|y\|^{\frac{1}{3}}_{[t,t+\pi]}.
\]
In other words, the function $\gamma$ in \eqref{eq:estAutxy} is given by $\gamma(s)=4s+s^{\frac{1}{3}}$.

We can now use the results in Proposition \ref{prop_nearest_control} to practically stabilize $\Sigma_{\text{ex}}$.
We choose the control set to be $\mathcal U_{\rm ex}$ given in \eqref{ex:U}, and the desired stability margin to be $\epsilon=1$. Then, based on the function $\gamma$ computed for the system $\Sigma_{\textrm{ex}}$, we get $\gamma(\delta) < 1$ if $\delta \in \left(\ 0,\  \frac{1}{8}\ \right]$. 
Using the same discrete set as in \eqref{ex:U} along with the function $\phi$ as in \eqref{phi_eq}, we can fix $\theta_{\rm ex}=0$ and $\alpha=0.1$ such that the system $\Sigma_{\rm ex}$ is globally practically stable with respect to $\mathbb B_\epsilon$, with $\epsilon = 1$, as shown in the simulation results in Figure~\ref{fig:ex1}. 

\end{example}

\section{Nearest-Neighbor Control for Incrementally Passive Systems With Constant Inputs}\label{sec:inc}

In many cases, the desired equilibrium point of the passive nonlinear system $\Sigma$ as in \eqref{systems_eq} is not equal to the minimum of the associated storage function $H$. Instead, it may correspond to an arbitrary constant input. For these cases, a constant input $u^*\in\R^m$ with its corresponding steady-state solution $x^*\in\R^n$ defines the steady-state relation given by the set
\begin{equation}
    \mathcal{E}:=\left\{(x^*,u^*)\in\R^n\times\R^m\, \bigg|\,  0=f(x^*)+g(x^*)u^*\right\}. \label{eq:eq_set}
\end{equation}

The problem of practically stabilizing the system $\Sigma$ around $x^*\in\R^n$ is equivalent to practically stabilizing $\overline{x} = x - x^*$ around the origin, with $\overline{(\cdot)}=(\cdot)-{(\cdot)}^*$ denoting the incremental variable. Thus, the incremental system is given by
\begin{equation}\label{eq:shifted_sys_aut}
    \overline{\Sigma}:
    \left\{
        \begin{array}{rl}
           \dot{\overline{x}}  & = \overline{f}(\overline{x}) + g(\overline{x}+x^*)\overline{u},\\
            \overline{y} & = h(\overline{x}+x^*)-h(x^*),
        \end{array}\right.
\end{equation}
with
$
\overline{f}(\overline{x}) = f(\overline{x}+x^*)-f(x^*) + \bracket{g(\overline{x}+x^*) -g(x^*)} u^*.
$ 
For this matter, the passivity of the mapping $\overline{u}\mapsto \overline{y}$ is, in the original system $\Sigma$, referred to as incremental passivity with respect to constant input; and is defined as follows \citep{Jayawardhana2007}.

\begin{defn}[Constant Incremental Passivity]
Consider the nonlinear system $\Sigma$ as in \eqref{systems_eq}. The system $\Sigma$ is said to be incrementally passive with respect to constant input if, for every $(x^*,u^*)\in\mathcal{E}$, the corresponding incremental system $\overline \Sigma$ in \eqref{eq:shifted_sys_aut} with input $\overline u$ and output $\overline y$, is passive; that is, there exists a storage function $H_0:\R^n\to\R_{\geq 0}$ such that
\begin{equation}
    \dot{H}_0 = \langle \nabla H_0, \dot{\overline x} \rangle \leq \langle{\overline{u}}, \overline{y}\rangle.
\end{equation}
\end{defn}

{\AT Note that the incremental passivity is a stronger requirement than the passivity notion considered in Section~\ref{sec:result}. In particular, one can find examples of systems which are passive but not incrementally passive.} Also,  constant incremental passivity defined above is equivalent to shifted passivity as in \citep{Monshizadeh2019,Schaft2016} and equilibrium-independent passivity as in \citep{Hines2011}. Nevertheless, the term constant incremental passivity is preferred in this paper because the pair $(x^*,u^*)$ can be arbitrary and most importantly, the incremental function is used in the definition. In the remainder of this section, we study stabilization of incrementally passive systems with finite set of control actions.

\subsection{Steady-state $u^* \in \cU$}

In the case of constant incremental passivity, the corresponding constant input $u^*$ is often known from the knowledge of the nominal system~\eqref{systems_eq}. Then we can simply design the finite input set $\mathcal{U}$ such that it contains $u^*$. Thus it is natural to adapt the assumption (A1) to the current setting that brings us to the following proposition. 

\begin{prop}\label{prop:phi_incrementalpassivity}
Consider the system $\Sigma$ as in \eqref{systems_eq}, and a finite set of control actions $\cU =\{u_0, u_1, \dots, u_p\}  \subset \R^m$.
Assume that:
\begin{enumerate}[leftmargin=3em]
    \item[(A2)] $\Sigma$ is constant-incrementally passive with the proper storage function $H_0(x,x^*)$ for all pair $(x^*,u^*)\in\mathcal{E}$;
    \item[(A3)] $u^*\in\cU$, with $u_0=u^*$, and there exists a subset $\mathcal{V}$ of $\mathcal{U}$ such that $u^*\in\rm{int}\bracket{\rm{conv}\bracket{\mathcal{V}}}$; and
    \item[(A4)] the autonomous incremental system $\overline{\Sigma}$ with $u=u^*$
    is large-time initial-state norm-observable, i.e.\ there exists $\tau>0$ and $\bar{\gamma}\in\mathcal{K}_\infty$ such that the solution of the autonomous incremental system $\overline{\Sigma}_{u=u^*}$ satisfies
    \[\|\overline{x}(t)\|\leq \bar{\gamma}\bracket{{\|\overline{y}\|}_{[t,t+\tau]}}\]
    for all $\overline{x}(0) \in \R^n, t\geq0$.
\end{enumerate}
Furthermore, for a given $\epsilon>0$, assume that $\bar{\gamma}(\delta)\leq \epsilon$, where $\delta > 0$ is the smallest number that satisfies
\begin{equation}\label{eq:defCvDelta}
V_{\cU} (u^*) \subseteq \mathbb{B}_\delta (u^*).
\end{equation}
Then the control law $u=\phi\bracket{\overline{y}-u^*}$, with 
\startmodif
$\phi:\R^m \rightrightarrows \cU$
\stopmodif
defined in \eqref{phi_eq}, globally practically stabilizes $\Sigma$ with respect to $\mathbb{B}_{\epsilon}\bracket{x^*}$.
\end{prop}

\startmodif
The proof of Proposition~\ref{prop:phi_incrementalpassivity} can be developed similarly to the proof of Proposition~\ref{prop_nearest_control}, by considering
\begin{equation}\label{eq:barPhieq}
    \overline{\phi}(\overline{y})= \phi(\overline{y}-u^*)-u^*
\end{equation}
with
\begin{equation}\label{eq:phi_bar}
   \overline{\phi}(\eta):=\argmin_{\overline{v}\in\overline{\mathcal{U}}}\BR{\norm{\overline{v}+\eta}}
\end{equation}
where the set
\begin{equation}\label{eq:U_bar}
    \overline{\mathcal{U}}:=\{\overline{v}\in\R^m\ \mid\ \overline{v}=v-u^*; v\in\mathcal{U}\}
\end{equation}
is defined by \emph{shifting} the original input set $\mathcal{U}$ such that $u^*$ is now the origin of the input/output space of the constant incremental system. This means that we can use the constant-incremental nearest-neighbor map $\overline{\phi}$ so that the constant incremental system has the same structure as \eqref{systems_eq}. Then the rest of the proof follows from the proof of Proposition~\ref{prop_nearest_control}. Finally, since the output and state variables of the constant incremental system converge to $\mathbb B_\delta$ and $\mathbb B_\epsilon$, respectively, as $t\to\infty$, we can conclude practical stability, i.e.\ $\overline y\to\mathbb B_\delta$ and $x\to\mathbb B_\epsilon(x^*)$ as $t\to\infty$.
\stopmodif

\subsection{Sector bounded feedback}

Similar to the results in the previous section, sector bounded nonlinear mapping $F$ that satisfies \eqref{eq:condSec} can easily be included in the constant-incrementally passive systems case. This is due to the fact given by \eqref{eq:barPhieq}. Then the following proposition is true.

\begin{prop}
Consider a nonlinear system $\Sigma$ described by \eqref{systems_eq} that satisfies (A2) and (A4); and a discrete set $\cU\subset\R^m$ satisfying (A3) so that \eqref{eq:defCvDelta} holds for some $\delta>0$. Let $\phi$ be as given in \eqref{phi_eq}; and let $\mu_{{\rm min},1}\in(0,1]$ be such that \eqref{theta1max_eq} holds for all $z\in\R^m$. Assume that \eqref{sector_condition_F} holds with the mapping $F$, along with constants $k_1,k_2,k_3$, satisfying \eqref{eq:condSec}. For a given $\epsilon>0$, assume that
\[
\bar{\gamma}\bracket{\delta/k_1}\le \epsilon.
\]
Then, the control law $u=\phi\bracket{F(\bar y)-u^*}$ globally practically stabilizes $\Sigma$ with respect to $\mathbb{B}_\epsilon(x^*)$.
\end{prop}

\subsection{Revisiting an illustrative example}

\begin{example}[continues=ex:1]\label{ex:2}
\normalfont
Consider the nonlinear system $\Sigma_{\text{ex}}$ along with the associated storage function $H(x)$ as in Example~\ref{ex:1}. 
It can be shown easily (following the main results in~\citep{Jayawardhana2007}) that $\Sigma_{\text{ex}}$ is constant-incrementally passive. Indeed, for any $(x^*,u^*)\in\mathcal{E}$, we can define
    \startmodif
\begin{align*}
    H_0(x,x^*)&=\frac{1}{2} \bracket{x_1-x_1^*}^2 + \frac{1}{2} \bracket{x_2-x_2^*}^2\\
    &\qquad\qquad+ \frac{1}{4} \bracket{x_3^4+3 {x_3^*}^4}-x_3 {x_3^*}^3, 
\end{align*}
which has a global unique minimum at $x^*$ and is related to the original storage function $H(x)$ following \citep{Jayawardhana2007} by $
H_0(x,x^*) = H(x) - H(x^*) - (x-x^*)^T\nabla H(x^*)$. 
    \stopmodif
It follows immediately that $\dot H_0 = \langle \overline{y},\overline{u}\rangle$.

We will now show that the autonomous incremental system of $\Sigma_\text{ex}$ satisfies the large-time initial-state norm observability conditions. Let the function $\overline{\gamma}$ be computed by considering $u=u^*$ for all $t\geq0$ as provided in the following. Consider the incremental system of $\Sigma_{\text{ex}}$ with $u=u^*$ for all $t\geq 0$, i.e.\

\startmodif
\begin{equation}
\overline{\Sigma}_{\text{ex}} \biggm|_{u=u^*} \ : \ \left\{    \begin{array}{rl}
\dot{\overline{x}}         & = \bbm{-x_2+x_2^*\\x_1-x_1^*\\-x_1+x_1^*} + \bbm{\overline{y}_2\\0\\0}, \\
\overline{y}         & = \bbm{x_1-x_1^*\\x_3^3-{x_3^*}^3}.
    \end{array} \right.
\end{equation}
\stopmodif

Following the computation in Example~\ref{ex:1}, we first compute the bound on the subsystem $\sbm{\overline x_1\\\overline x_2}$ by considering
\startmodif
$\overline{y}_2$ as the input and $\overline{y}_1 = \overline{x}_1$ as the output.
\stopmodif
Then we have a linear system with $\overline{A}=A=\sbm{0 & -1\\1 & 0}$, $\overline{B}=B=\sbm{1\\0}$, $\overline{C}=C=\bbm{1&0}$. Hence, following a similar routine computation as before, we get

\begin{align*}
    \left\|\sbm{\overline x_1(t)\\\overline x_2(t)}\right\| & \leq 2 \left(\|\overline y_1\|_{[t,t+\pi]} + \|\overline y_2\|_{[t,t+\pi]} \right) 
    \\  & \leq 4 \|\overline y\|_{[t,t+\pi]}.
\end{align*}

Accordingly, for $\overline{x}_3$, we have that
$
\overline x_3 = \frac{\overline{y}_2}{x_3^2+{x_3^*}^2+x_3 x_3^*}.
$ 
For 
\startmodif
any 
\stopmodif
$x_3^*\neq 0$, we have that $x_3^2+{x_3^*}^2+x_3 x_3^*\geq \frac{3}{4}{x_3^*}^2$, for all $x_3$. Hence,
    \startmodif
\begin{align*}
    \|\overline{x}(t)\| \leq \left\|\sbm{\overline x_1(t)\\\overline x_2(t)}\right\| + \|\overline x_3\|
     & \leq 4 \|\overline y\|_{[t,t+\pi]} + \frac{4}{3{x_3^*}^2}{\|\overline y_2\|}_{[t,t+\pi]}\\
    &\leq 4 \|\overline y\|_{[t,t+\pi]} + \frac{4}{3{x_3^*}^2}{\|\overline y\|}_{[t,t+\pi]}.
\end{align*}
    \stopmodif

In other words, the large-time initial-state norm-observability function is given by $\overline{\gamma}(s,{x_3^*}) = 4s+\frac{4}{3{x_3^*}^2}s^2$.

We can now use the results in Proposition \ref{prop:phi_incrementalpassivity} to practically stabilize $\Sigma_{\text{ex}}$ around any arbitrary steady-state relation $(x^*,u^*)\in\mathcal{E}$. Fix $x^*={\sbm{0 & 0 & -1}}^\top$, $u^*={\sbm{1 & 0}}^\top$, and $\epsilon=0.5$. Then, by the large-time initial-state norm-observability property of the incremental system, we can choose $\delta=0.1$ to generate the discrete set of control actions. In this case, we can translate the previously used discrete set such that $u^*$ is in the realizable control actions, i.e. $\overline{\mathcal{U}}_{ex} := \mathcal{U}_{ex} + u^*$ with $\mathcal{U}_{ex}$ be as the discrete input set used in Example~\ref{ex:1}. The mapping control law with the mapping $\phi$ can then be demonstrated as shown in Figure~\ref{fig:ex2}. 

\begin{figure}
    \centering
    \includegraphics[width=0.95\linewidth]{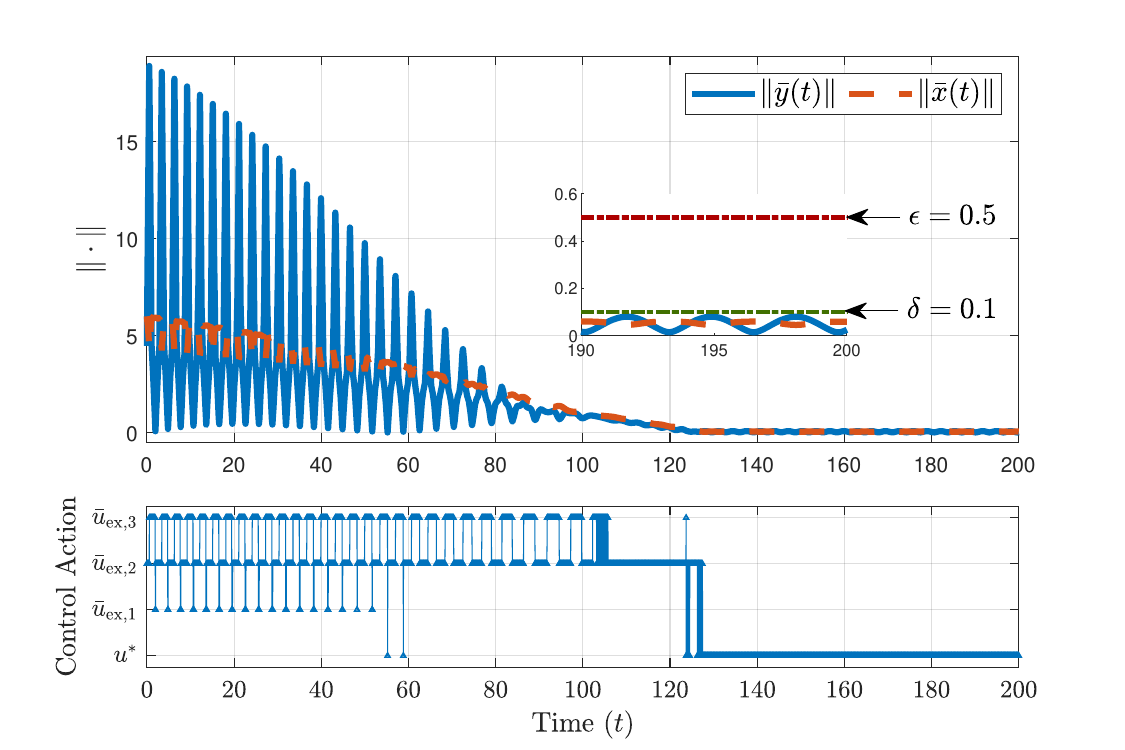}
    \caption{\small Simulation results of $\Sigma_{ex}$ using the control approach proposed in the Propostion~\ref{prop:phi_incrementalpassivity} with discrete input set $\overline{\mathcal{U}}_{ex} := \mathcal{U}_{ex} + u^*$ with $\mathcal{U}_{ex}$ be as the input set used in Example~\ref{ex:1}. Here we have that $u^*\in\overline{\mathcal{U}}_{ex}$. Similar to before, in this simulation, once both the state $x$ and the output $y$ enter their respective convergence ball, the control action is switched to $u^*$ for the rest of the simulation.}
    \label{fig:ex2}
\end{figure}

\end{example}

\section{Minimal Control Actions: Constructions and Bounds}\label{sec:minimal}

In the earlier sections, we have shown that a nearest neighbor approach is a powerful tool for global practical stabilization of passive nonlinear systems. Indeed, when we are given a limited choice of static control inputs, assumptions (A1) and (A3) provide us a way to check whether the given set of inputs is applicable by means of nearest neighbor approach for the practical stabilization problem. If these assumptions hold for a finite set $\cU$, then it is of interest to compute the smallest number $\delta >0$ associated with Voronoi cell $V_{\mathcal{U}}(u^*)$, such that
\startmodif
\[
V_{\cU}(u^*) \subset \mathbb{B}_\delta(u^*).
\]
\stopmodif
Since our control design achieves convergence up to a ball of radius $\gamma(\delta)$, with $\gamma(\cdot)$ being the output-to-state gain in large-time initial-state norm-observability assumption, the knowledge of $\delta$ basically determines how close the trajectories can get to the desired equilibrium with our proposed controller.
Let us recall the assumption (A1) and its generalization (A3), where we assume that, for a finite set $\cU$, the desired equilibrium $u^* \in {\rm int}({\rm conv}(\cV))$. To obtain $\cU$ of minimal cardinality, the following result, borrowed from \startmodif\cite[Corollary~9.5]{Brondsted1983}\stopmodif, is of interest:
\begin{lem}\label{thm:minimal_card}
For a finite set $\mathcal{S}\subset\R^m$, the minimal cardinality of $\mathcal{S}$ such that $\rm{int}\bracket{\rm{conv}\bracket{\mathcal{S}}}\neq \emptyset$ is equal to $m+1$.
\end{lem}

An immediate consequence is that, for practical stabilization of passive systems, it suffices to consider a control set $\cU$ with 
\startmodif
$m+2$ elements (including $u^*$), 
\stopmodif
provided they satisfy a certain geometric configuration.

\begin{cor}\label{cor:minimal_NNC_stabilizing_input}
Let the set $\mathcal{U}$ be such that $u^*\in{\rm int}\bracket{{\rm conv}\bracket{\mathcal{U}}}$. If ${\rm conv}\bracket{\mathcal{U}\setminus\{u^*\}}$ is an $m$-simplex, then $\mathcal{U}$ is a minimal set that satisfies (A3).
\end{cor}

In the remainder of this section, we will work with two particular choices of the set $\cU$ with cardinality 
\startmodif
$m+2$
\stopmodif
that satisfy (A1) or (A3). We give a closed-form expression of $\delta$ for these sets in terms of the elements $\cU$. For the sake of simplicity, we fix $u^* = 0$ in these computations. The two cases we consider are:
\begin{enumerate}
\item The set $\cU = \cS_{\rm reg} \cup \{0\}$, where
\begin{equation}\label{eq:regSimp}
\mathcal{S}_{\rm{reg}}=\lambda\BR{e_1,\dots,e_m,\frac{1-\sqrt{m+1}}{m}\mathbbm{1}},
\end{equation}
for some $\lambda > 0$, and the barycenter of $\mathcal{S}_{m,\rm{reg}} = {\rm conv} (\cS_{\rm reg})$ is
\begin{equation}\label{eq:baryCent}
b_{\mathcal{S}_{\rm reg}} = \lambda \frac{\sqrt{m+1}-1}{m\sqrt{m+1}}\mathbbm{1}.
\end{equation}
\item We take $\cU = \cS_{\rm reg}^0 \cup \{ 0 \}$, where $\cS_{\rm reg}^0 = \cS_{\rm reg} - b_{\mathcal{S}_{\rm reg}}$ and ${\rm conv} (\cS_{\rm reg}^0)$ has barycenter at the origin.
\end{enumerate}

In the next two lemmas, we basically compute a bound on the sets 
\startmodif
$V_{\cS_{\rm reg} \cup \{ 0 \}} (0)$ and $V_{\cS_{\rm reg}^0 \cup \{ 0 \}} (0)$. 
\stopmodif
It is noted that the results apply to the case when $u^* \neq 0$ since the set
\startmodif
$\cV=(\mathcal{S}_{\rm{reg}}\cup\{0\})+u^*$ (or $\cV=(\mathcal{S}_{\rm{reg}}^0\cup\{0\})+u^*$)
\stopmodif
is such that $u^*\in{\rm int}\bracket{{\rm conv}\bracket{\cV}}$ and the relative position of $u^*$ with respect to the vertices of $\cV$ is the same as the relative position of the origin with respect to the vertices of $\mathcal{S}_{\rm{reg}}$ (or $\mathcal{S}_{\rm{reg}}^0$); hence it has the same bound.

\begin{lem}\label{lemma:bound_sreg}
    For the set $V_{\mathcal{S}_{\rm{reg}}\cup\{0\}}(0)$, the smallest bound $\delta>0$ such that $V_{\mathcal{S}_{\rm{reg}}\cup\{0\}}(0)\subset\mathbb{B}_\delta$ is given by 
    \startmodif
    \[\delta=\left\{\begin{array}{cl}
        \frac{\lambda }{2}, & \text{if }m=1, \\
        \frac{\lambda }{2} \sqrt{m-1+{(2-m-\sqrt{m+1})}^2},& \text{otherwise.}
    \end{array}\right.
    \]
    \stopmodif
\end{lem}

\begin{proof}
First, we observe that the set $V_{\mathcal{S}_{\rm{reg}}\cup\{0\}}(0)$ is equal to the solution set of the following inequalities. The vector $x=\bbm{x_1 & \dots & x_m}^\top \in V_{\mathcal{S}_{\rm{reg}}\cup\{0\}}(0)$ if 
\begin{align}
    x_i&\leq \frac{\lambda }{2},\quad\quad\quad i=1,\dots,m,\label{ineq:1_lem_reg_simplex}\\
    \frac{1-\sqrt{m+1}}{m}\mathbbm{1}^\top x&\leq \lambda\frac{\bracket{1-\sqrt{m+1}}^2}{2m}. \label{ineq:2_lem_reg_simplex}
\end{align}
Since $m\geq 1$, the inequality \eqref{ineq:2_lem_reg_simplex} can simply be rewritten as
\[\mathbbm{1}^\top x\geq \frac{\lambda }{2}\bracket{1-\sqrt{m+1}}.\]
{\color{black}Next, we observe that each of the vertices of the Voronoi cell $V_{\mathcal{S}_{\rm{reg}}\cup\{0\}}(0)$ can be obtained by solving $m$ equations taken from \eqref{ineq:1_lem_reg_simplex} and/or \eqref{ineq:2_lem_reg_simplex}.} Let $\mathcal{V}$ be the set of all vertices of $V_{\mathcal{S}_{\rm{reg}}\cup\{0\}}(0)$. Then $\mathcal{V}=\{\frac{\lambda }{2} \mathbbm{1}\}\bigcup\limits_{i=1}^m\{\frac{\lambda }{2} \tilde{v}_i\}$ with $\tilde{v}_i$ being a column vector where the $i$-th element is given by $2-m-\sqrt{m+1}$ and the other $m-1$ elements are $1$. Therefore, the minimum value of $\delta$ for which $V_{\mathcal{S}_{\rm{reg}}\cup\{0\}}(0) \subset \mathbb{B}_\delta$ is given by
\startmodif
\[\delta_{m=1}=\maxxx\limits_{\tilde{v}\in\mathcal{V}}\{ \|\tilde{v}\| \}=\frac{\lambda }{2}\|1\|=\frac{\lambda }{2}\]
and
\[\begin{split}
    \delta_{m>1}&=\maxxx\limits_{\tilde{v}\in\mathcal{V}}\{ \|\tilde{v}\| \} = \frac{\lambda }{2}\|\tilde v_i\|\\&=\frac{\lambda }{2}\sqrt{m-1+{(2-m-\sqrt{m+1})}^2}.
\end{split}
\]
\stopmodif
which is the desired expression.
\end{proof}

Next, let us consider the regular $m$-simplex centered at the origin with vertices $\mathcal{S}_{{\rm reg}}^0$.

\begin{lem}\label{lemma:bound_sreg_0}
Consider the vertices of a regular $m$-simplex centered at the origin $\mathcal{S}_{\rm{reg}}^0 = \cS_{\rm reg} - b_{\cS_{\rm reg}}$ where $\cS_{\rm reg}$ and $b_{\cS_{\rm reg}}$ are defined in \eqref{eq:regSimp} and \eqref{eq:baryCent}.
Then the bound $\delta>0$ such that $V_{\mathcal{S}_{\rm{reg}}^0\cup\{0\}}(0)\subset\mathbb{B}_\delta$ is given by
\[
\delta=\lambda \frac{m}{2} \sqrt{\frac{m}{m+1}}.
\]
\end{lem}

\begin{proof}
Let us denote the set $\mathcal{S}:=\mathcal{S}_{\rm{reg}}^0\bigcup\{0\}$. Then, by following the same proof as before, we have that the set $V_{\mathcal{S}}(0)$ is equal to the solution set of the following system of inequalities,
\startmodif
\begin{align}
    \bracket{e_i-\frac{\sqrt{m+1}-1}{m\sqrt{m+1}}\mathbbm{1}}^\top x&\leq \frac{\lambda }{2}\frac{m}{m+1},\quad i=1,\dots,m,\label{ineq:3_lem_reg_simplex}\\
    -\frac{1}{\sqrt{m+1}}\mathbbm{1}^\top x&\leq \frac{\lambda }{2}\frac{m}{m+1}.\label{ineq:4_lem_reg_simplex}
\end{align}
\stopmodif
{\color{black}Since all points in $\mathcal{S}_{\rm{reg}}^0$ have the same distance from the origin, we can pick any set of $m$ equations from the above inequalities in order to get one of the vertices of $V_{\mathcal{S}}(0)$. Let us now choose all $m$ equations from \eqref{ineq:3_lem_reg_simplex} because they have a nice symmetric structure given by,}
\begin{equation}
    \underbrace{\left[\begin{smallmatrix}a & -1 & \cdots & -1\\
    -1 & a & \cdots & -1\\
    \vdots & \vdots & \ddots & \vdots\\
    -1 & -1 & \cdots & a\end{smallmatrix}\right]}_\cA v=\frac{\lambda}{2}\frac{m^2}{m+1-\sqrt{m+1}}\mathbbm{1},
\end{equation}
where $a=m+\sqrt{m+1}$. Since $\cA$ is symmetric, we can find a symmetric matrix $\cA^{-1}$ such that $\cA\cA^{-1}=\cA^{-1}\cA=I_{m\times m}$. Via routine computation, we have that, {\color{black} $\cA^{-1} = bI + c (\mathbbm{1}\mathbbm{1}^\top - I)$, 
whose main diagonal elements are
$
b=\frac{\sqrt{m+1}(m-2)+2}{m^2\sqrt{m+1}},
$ 
and its off-diagonal elements are
$
c=\frac{1}{2(1+\sqrt{m+1})+m(2+\sqrt{m+1})}.
$
}
Thus, the point $v$ is given by
\[v=\frac{\lambda}{2}\frac{m}{\sqrt{m+1}}\mathbbm{1}.\]
Therefore, the minimum bound on the set $V_{\mathcal{S}}(0)$ is,
\begin{align*}
    \delta = \|v\| = \frac{\lambda}{2}\frac{m}{\sqrt{m+1}}\|\mathbbm{1}\| 
    =\lambda \frac{m}{2} \sqrt{\frac{m}{m+1}}.
\end{align*}
which completes the proof.
\end{proof}

We have shown in Lemma~\ref{lemma:bound_sreg} and Lemma~\ref{lemma:bound_sreg_0} above that for the two types of discrete sets, whose elements form the vertices of  regular $m$-simplices', the minimum bounds of the Voronoi cell of the origin can be computed in a closed-form manner. Now, for a given incrementally passive system $\Sigma$ and admissible reference signal $u^*$ with large-time norm-observability function $\gamma\in\mathcal{K}$ when $u=u^*$, for a given stability margin $\epsilon>0$, the value of the bound $\delta$ can be chosen as large as possible such that $\gamma\bracket{\delta}\leq \epsilon$. Thus, for a given $\epsilon>0$, norm-observability function $\gamma$ of the system $\Sigma$, and a desired rotation matrix $R$, we can choose $\delta>0$ that satisfies $\gamma\bracket{\delta}\leq \epsilon$ and construct the minimal set $\mathcal{U}\subset\R^m$ that satisfies (A3) as follows:
\begin{enumerate}
    \item $\mathcal{U}:=\bracket{R\mathcal{S}_{\rm{reg}}\cup\{0\}}+u^*$ with
\startmodif
    \[\lambda=\min\left\{2 \delta,\frac{2 \delta}{\sqrt{m-1+\bracket{2-m-\sqrt{m+1}}^2}}\right\},\]
    \stopmodif
    or;
    \item \startmodif $\mathcal{U}:=\bracket{R\mathcal{S}^0_{\rm{reg}}\cup\{0\}}+u^*$ 
    \stopmodif
    with
    \[\lambda= \frac{2\delta}{m}\sqrt{\frac{m+1}{m}}.\]
\end{enumerate}

\begin{example}
Recall the discrete set $\mathcal{U}_{\rm ex}$ as in Example~\ref{example:exU}. The same discrete set can be constructed by using $\mathcal{U}:=\bracket{R\mathcal{S}_{\rm{reg}}(0)\cup\{0\}}+u^*$; by fixing $\alpha=\delta$ and 
\[R=-\frac{\sqrt{2}}{2}\bbm{\ \ \ \sin{\theta_{\rm ex}}+\cos{\theta_{\rm ex}}\ \ & \sin{\theta_{\rm ex}}-\cos{\theta_{\rm ex}}\\
-\sin{\theta_{\rm ex}}+\cos{\theta_{\rm ex}}\ \ & \sin{\theta_{\rm ex}}+\cos{\theta_{\rm ex}}}.\]
\end{example}

\section{Conclusions and Further Research}\label{sec:conc}

We have considered practical stabilization of continuous-time (incrementally) passive nonlinear systems using output-feedback where the control inputs only take values among the available actions in a given finite discrete set. 
We propose simple ways to select the control actions at each time instance where we have shown that our proposed control laws are able to stabilize the systems up to some desirable distance from the equilibrium. In addition, our results provide an insight on the lower bound on the number of control elements that guarantee practical stability. We have also provided 
methods to design the finite set of control actions with minimal cardinality. 
Questions related to improving the convergence rate with more (than necessary) control elements and/or to eliminate the chattering effects are being investigated as further directions of research.

\section{Acknowledgements}
{\AT The authors would like to thank anonymous reviewers for their helpful and constructive comments.}

The project of Muhammad Zaki Almuzakki is fully supported by \emph{Lembaga Pengelola Dana Pendidikan Republik Indonesia} (LPDP-RI) under contract No. PRJ-851/LPDP.3/2016.

\bibliography{bib_all}         

\end{document}